\documentclass{elsart}

\journal{}

\usepackage{amsmath,amssymb}

\newcommand{\ketsor}[2]{\genfrac{}{}{0pt}{2}{#1}{#2}}

\begin{document}

\newcommand {\stirlingf}[2]{\genfrac[]{0pt}{}{#1}{#2}}
\newcommand {\stirlings}[2]{\genfrac\{\}{0pt}{}{#1}{#2}}

\newtheorem{Theorem}{Theorem}

\begin{frontmatter}

 \title{On the coefficients of power sums of arithmetic progressions}

 \author[A]{Andr\'as Bazs\'o} and
 \ead{bazsoa@science.unideb.hu}
 \author[I]{Istv\'an Mez\H{o}}
 \ead{istvanmezo81@gmail.com}

 \address[A]{Institute of Mathematics, MTA-DE Research Group "Equations, Functions and Curves", Hungarian Academy of Sciences and University of Debrecen, P.O. Box 12, H-4010 Debrecen, Hungary}
 \address[I]{Department of Mathematics, Nanjing University of Information Science and Technology, Nanjing, 210044 China}

\begin{abstract}
We investigate the coefficients of the polynomial
\[ S_{m,r}^n(\ell)=r^n+(m+r)^n+(2m+r)^n+\cdots+((\ell-1)m+r)^n. \] We prove that these can be given in terms of Stirling numbers of the first kind and $r$-Whitney numbers of the second kind. Moreover, we prove a necessary and sufficient condition for the integrity of these coefficients. 
\end{abstract}

\begin{keyword}arithmetic progressions, power sums, Stirling numbers, $r$-Whitney numbers, Bernoulli polynomials
\MSC{11B25, 11B68, 11B73}
\end{keyword}

\end{frontmatter}

\maketitle

\section{Introduction}

Let $n$ be a positive integer, and let
\[S_n(\ell)=1^n+2^n+\cdots+(\ell-1)^n\]
be the power sum of the first $\ell-1$ positive integers. It is well known that $S_n(\ell)$ is strongly related to the Bernoulli polynomials $B_n (x)$ in the following way
\[S_n(\ell)=\frac{1}{n+1}(B_{n+1}(\ell)-B_{n+1}).\]
where the polynomials $B_n(x)$ are defined by the generating series 
\[\frac{t e^{tx}}{e^t-1}=\sum_{k=0}^{\infty}B_{k}(x)\frac{t^k}{k!}\]
and $B_{n} = B_{n} (0)$ is the $n$th Bernoulli number. 

It is possible to find the explicit coefficients of $\ell$ in $S_n(\ell)$ \cite{GKP}:
\begin{equation}
S_n(\ell)=\sum_{i=0}^{n+1}\ell^i\left(\sum_{k=0}^nS_2(n,k)S_1(k+1,i)\frac{1}{k+1}\right),\label{Spolcoeff}
\end{equation}
where $S_1(n,k)$ and $S_2(n,k)$ are the (signed) Stirling numbers of the first and second kind, respectively.

Recently, Bazs\'o et al. \cite{BPS} considered the more general power sum
\[S_{m,r}^n(\ell)=r^n+(m+r)^n+(2m+r)^n+\cdots+((\ell-1)m+r)^n,\]
where $m\neq 0, r$ are coprime integers. Obviously, $S_{1,0}^n(\ell) = S_n (\ell)$.
They, among other things, proved that $S_{m,r}^n(\ell)$ is a polynomial of $\ell$ with the explicit expression
\begin{equation} \label{eq:Smrn}
S_{m,r}^n(\ell)=\frac{m^n}{n+1}\left(B_{n+1}\left(\ell+\frac{r}{m}\right)-B_{n+1}\left(\frac{r}{m}\right)\right).
\end{equation}
In \cite{How}, using a different approach, Howard also obtained the above relation via generating functions. Hirschhorn \cite{Hirschhorn} and Chapman \cite{Chapman} deduced a longer expression which contains already just binomial coefficients and Bernoulli numbers.

For some related diophantine results on $S_{m,r}^n(\ell)$ see \cite{BBKPT,GyPsurv,R,R2,BKLP} and the references given there.

Our goal is to give the explicit form of the coefficients of the polynomial $S_{m,r}^n(\ell)$, thus generalizing \eqref{Spolcoeff}. In this expression the Stirling numbers of the first kind also will appear, but, in place of the Stirling numbers of the second kind a more general class of numbers arises, the so-called $r$-Whitney numbers introduced by the second author \cite{Mezo}.

The $r$-Whitney numbers $W_{m,r}(n,k)$ of the second kind are generalizations of the usual Stirling numbers of the second kind with the exponential generating function
\[\sum_{n=k}^\infty W_{m,r}(n,k)\frac{z^n}{n!}=\frac{e^{rz}}{k!}\left(\frac{e^{mz}-1}{m}\right)^k.\]

For algebraic, combinatoric and analytic properties of these numbers see \cite{CJ,Rahmani} and \cite{CC1,CC2}, respectively.

First, we prove the following.

\begin{Theorem} For all parameters $\ell>1,n,m>0,r\ge0$ we have
\[S_{m,r}^n(\ell)=\sum_{i=0}^{n+1}\ell^i\left(\sum_{k=0}^n\frac{m^kW_{m,r}(n,k)}{k+1}S_1(k+1,i)\right).\]
\end{Theorem}

\textit{Proof.}\ 
The formula which connects the power sums and the $r$-Whitney numbers is the next one from \cite{Mezo}:
\[(mx+r)^n=\sum_{k=0}^nm^kW_{m,r}(n,k)x^{\underline{k}}.\]
Here $x^{\underline{k}}=x(x-1)\cdots(x-k+1)$ is the falling factorial. We can see that it is enough to sum from $x=0,1,\dots,\ell-1$ to get back $S_{m,r}^n(\ell)$. Hence
\[S_{m,r}^n(\ell)=\sum_{k=0}^nm^kW_{m,r}(n,k)\sum_{x=0}^{\ell-1}x^{\underline{k}}.\]
The inner sum can be determined easily (see \cite{GKP}):
\[\sum_{x=0}^{\ell-1}x^{\underline{k}}=\frac{\ell^{\underline{k+1}}}{k+1}+\delta_{k,0}.\]
The Kronecker delta will never appear, because if $k=0$ then the $r$-Whitney number is zero (unless the trivial case $n=0$, which we excluded). Therefore, as an intermediate formula, we now have that
\[S_{m,r}^n(\ell)=\sum_{k=0}^nm^kW_{m,r}(n,k)\frac{\ell^{\underline{k+1}}}{k+1}.\]
The falling factorial $\ell^{\underline{k+1}}$ is a polynomial of $\ell$ with Stirling number coefficients:
\[\ell^{\underline{k+1}}=\sum_{i=0}^{k+1}S_1(k+1,i)\ell^i.\]
Substituting this to the formula above, we obtain:
\[S_{m,r}^n(\ell)=\sum_{k=0}^n\frac{m^kW_{m,r}(n,k)}{k+1}\sum_{i=0}^{k+1}S_1(k+1,i)\ell^i.\]
Since $S_1(k+1,i)$ is zero if $i>k+1$, we can run the inner summation up to $n+1$ (this is taken when $k=n$) to make the inner sum independent of $k$. Altogether, we have that
\[S_{m,r}^n(\ell)=\sum_{i=0}^{n+1}\ell^i\sum_{k=0}^n\frac{m^kW_{m,r}(n,k)}{k+1}S_1(k+1,i).\]
This is exactly the formula that we wanted to prove.
\hfill$\Box$

Now we give some elementary consequences of the theorem. The proofs are trivial.

\noindent \textbf{Remark.} The next properties of the polynomial $S_{m,r}^n(\ell)$ hold true for all parameters $\ell>1,n>0,r,m\ge0$:
\begin{enumerate}
	\item[(i)] The constant term of $S_{m,r}^n(\ell)$ is 0,
	\item[(ii)] The leading coefficient of $S_{m,r}^n(\ell)$ is $m^n/(n+1)$,
	\item[(iii)] $S_{m,r}^n(\ell)$ is a polynomial of $\ell$ of degree $n+1$ unless $m=0$; in this latter case the degree is $n$.
\end{enumerate}

The above statements also follow from \eqref{eq:Smrn}.

\section{The integer property of the coefficients in $S_{m,r}^n(\ell)$}

The coefficients of the polynomial $S_{m,r}^n(\ell)$ are not integer in the overwhelming majority of the cases:
\begin{align*}
S_{1,0}^1(\ell)=&\frac{\ell(\ell-1)}{2},\\
S_{2,5}^2(\ell)=&\frac13\ell(47+24\ell+4\ell^2),
\end{align*}
\begin{center}
etc.
\end{center}

However, we revealed that in special cases the polynomial $S_{m,r}^n(\ell)$ has integer coefficients. Several parameters are in the next table.

\begin{center}
\begin{tabular}{c|c|c}
$m$&$r$&$n$\\\hline\hline
2&1&3\\
2&3&3\\
2&5&3\\
4&3&3\\
4&5&3
\end{tabular}
\end{center}

For example,
\[S_{2,1}^3(\ell)=\ell^2(2\ell^2-1),\]
or
\[S_{2,3}^3(\ell)=\ell(2+\ell)(2\ell^2+4\ell+3).\]

From the formula of Theorem 1 it can be seen that if
\[(k+1)\mid m^kW_{m,r}(n,k)\quad (k=1,2,\dots,n),\]
then we get integer coefficients. 

To find another condition which is necessary and sufficient for the integrity of the coefficients in $S_{m,r}^n(\ell)$, we recall the following well known properties of Bernoulli polynomials and Bernoulli numbers.
\begin{equation} \label{eq:ber1}
B_n (x+y) = \sum_{k=0}^n{\binom{n}{k} B_k (x) y^{n-k}} = \sum_{k=0}^n{\binom{n}{k} B_k (y) x^{n-k}};
\end{equation}
\begin{equation} \label{eq:ber2}
B_n (x) =  \sum_{k=0}^n{\binom{n}{k} B_k x^{n-k}};
\end{equation} 
\begin{equation} \label{eq:ber3}
B_3 = B_5 = B_7 = \ldots = 0.
\end{equation}
By the \textit{denominator} of a rational number $q$ we mean the smallest positive integer $d$ such that $dq$ is an integer. We recall also the von Staudt theorem
\begin{equation} \label{eq:ber4}
\Lambda_{2n} = \prod_{\ketsor{(p-1) \mid 2n}{p \text{ prime}}} p,
\end{equation}
where $\Lambda_n$ is the denominator of $B_n$. In particular, $\Lambda_n$ is a square-free integer, divisible by $6$.
For the proofs of \eqref{eq:ber1}-\eqref{eq:ber3} see e.g. the work of Brillhart \cite{Brill}.

Let $2 \leq j \leq n$ be an even number and put
\begin{multline} \label{eq:fkj}
f(n,j) := \text{lcm} \left(\frac{\Lambda_j}{\gcd \left(\Lambda_j, \binom{n+1}{j}\binom{j}{j}\right)}, \frac{\Lambda_j}{\gcd \left(\Lambda_j, \binom{n+1}{j+1}\binom{j+1}{j}\right)}, \ldots, \right.\\
 \left. \frac{\Lambda_j}{\gcd \left(\Lambda_j, \binom{n+1}{n}\binom{n}{j}\right)} \right).
\end{multline}
Further, we define
\begin{equation} \label{eq:Fk}
F(n):= \begin{cases}
\text{lcm} \left(\text{rad}(n+1), f(n,2), f(n,4), \ldots, f(n,n)\right) & \  \text{ if $n$ is even,} \\
\text{lcm} \left(\text{rad}(n+1), f(n,2), f(n,4), \ldots, f(n,n-1)\right) & \  \text{ if $n$ is odd,} 
\end{cases}
\end{equation}
where 
\[ 
\text{rad}(n) = \prod_{\ketsor{p \mid n}{p \text{ prime}}} p .
\]

\begin{Theorem}
The polynomial $S_{m,r}^n(\ell)$ has integer coefficients if and only if $F(n) \mid m$.
\end{Theorem}

\textit{Proof.}\ 
By relations \eqref{eq:Smrn}, \eqref{eq:ber1} and \eqref{eq:ber2} we can rewrite $S_{m,r}^n(\ell)$ as follows:
\begin{align} 
S_{m,r}^n(\ell) =&  \frac{m^n}{n+1}\left(B_{n+1}\left(\ell+\frac{r}{m}\right)-B_{n+1}\left(\frac{r}{m}\right)\right) = \label{eq:prf1}\\
=& \frac{m^n}{n+1} \left(\left(\sum_{k=0}^{n+1}{\binom{n+1}{k} B_k \left(\frac{r}{m}\right) \ell^{n+1-k}}\right) - B_{n+1} \left(\frac{r}{m}\right)\right) = \label{eq:prf2}\\
=& \frac{m^n}{n+1} \sum_{k=0}^{n}{\binom{n+1}{k} B_k \left(\frac{r}{m}\right) \ell^{n+1-k}} = \label{eq:prf3}\\
=& \frac{m^n}{n+1} \sum_{k=0}^{n}{\binom{n+1}{k} \left(\sum_{j=0}^k \binom{k}{j} B_j \cdot \left(\frac{r}{m}\right)^{k-j}\right) \ell^{n+1-k}} \label{eq:prf4}
\end{align}

We denote the common denominator of the coefficients of $S_{m,r}^n(\ell)$ by $Q$. One can see from \eqref{eq:prf1} that the polynomial has integral coefficients if and only if $m$ is divisible by $Q$. Thus we have to determine $Q$. 

By \eqref{eq:prf4} we observe that neither $m$ nor $r$ occurs in $Q$. Moreover, the only algebraic expressions which may affect $Q$ in \eqref{eq:prf4} are on one hand $n+1$ and on the other hand, the denominators of the Bernoulli numbers involved, which are $2, \Lambda_j (2 \leq j \leq n \text{ even})$ by \eqref{eq:ber3} and the von Staudt theorem.    

It can easily be seen that $n+1 \mid m^n$ if $\text{rad}(n+1) \mid m$. Indeed, supposing the contrary, i.e., that $\text{rad}(n+1) \mid m$ and $n+1 \nmid m^n$, it implies that there is a prime factor $p$ of $n+1$ such that $p^{n+1}$ divides $n+1$. Hence $2^{n+1} \leq p^{n+1} \leq n+1$, which is a contradiction.

Let $2 \leq j \leq n$ be an even index. It follows from \eqref{eq:prf4} that the contribution of $\Lambda_j$ to the common denominator $Q$ is precisely $f(n,j)$ defined in \eqref{eq:fkj}. In other words, if $f(n,j) \mid m$, then every term of \eqref{eq:prf4} containing the factor $B_j$ has integer coefficients.

In conclusion, we obtained that $Q$ is the least common multiple of $\text{rad}(n+1)$ and $f(n,j)$ for all even $j \in \left[2,n\right]$, which number we denoted in \eqref{eq:Fk} by $F(n)$. The theorem is proved.
\hfill$\Box$

\noindent \textbf{Remark.} An easy consequence of our Theorem 2 is that $S_n(\ell)=S_{1,0}^n(\ell) \notin \mathbb{Z}[x]$ for any positive integer $n$.

Some small values of $F(n)$ are listed in the following table. These are results of an easy computation in MAPLE.

\begin{center}
\begin{tabular}{c|c||c|c||c|c||c|c}
$n$ & $F(n)$ & $n$ & $F(n)$ & $n$ & $F(n)$ & $n$ & $F(n)$ \\ \hline \hline
1 & 2 & 6 & 42 & 11 & 6 & 16 & 510 \\
2 & 6 & 7 & 6 & 12 & 2730 & 17 & 30 \\
3 & 2 & 8 & 30 & 13 & 210 & 18 & 3990 \\
4 & 30 & 9 & 10 & 14 & 30 & 19 & 210 \\
5 & 6 & 10 & 66 & 15 & 6 & 20 & 2310
\end{tabular}
\end{center}

\section*{Acknowledgements}

The first author was supported by the Hungarian Academy of Sciences and by the OTKA grant NK104208.

The research of Istv\'an Mez\H{o} was supported by the Scientific Research Foundation of Nanjing University of Information Science \& Technology, and The Startup Foundation for Introducing Talent of NUIST. Project no.: S8113062001


\begin{thebibliography}{AAAAAA}

\bibitem{BPS}
A. Bazs\'o, A. Pint\'er and H. M. Srivastava, On a refinement of Faulhaber's theorem concerning sums of powers of natural numbers, Appl. Math. Lett. 25 (2012), 486--489.

\bibitem{BKLP}
{ A.~Bazs\'o}, { D.~Kreso}, { F.~Luca} and {{\'A.}~Pint\'er},
  On equal values of power sums of arithmetic progressions, {Glas. Mat. Ser. III\/}, {47} (2012), 253--263. 

\bibitem{BBKPT}
{Y.~F. Bilu}, {B.~Brindza}, {P. ~Kirschenhofer}, {{\'A.}~Pint\'er} and {R. F. ~Tichy}, Diophantine equations and Bernoulli polynomials (with an Appendix by A. Schinzel), {Compositio Math.\/}, {131} (2002), 173--188.

\bibitem{Brill}
{J. ~Brillhart}, On the Euler and Bernoulli polynomials, { J. Reine Angew. Math.\/}, { 234} (1969),  45--64.

\bibitem{CJ}
G.-S. Cheon and J.-H. Jung, $r$-Whitney numbers of Dowling lattices, Discrete Math. 312(15) (2012),  2337--2348.

\bibitem{CC1}
R. B. Corcino and C. B. Corcino, On the maximum of generalized Stirling numbers, Util. Math. 86. (2011), 241-256.

\bibitem{CC2}
C. B. Corcino and R. B. Corcino, Asymptotic estimates for second kind generalized Stirling numbers, J. Appl. Math. (to appear).

\bibitem{Chapman}
R. Chapman, Evaluating $\sum_{n=1}^N(a+nd)^p$ again, Math. Gaz. {92} (2008), 92--94.

\bibitem{GKP}R. L. Graham, D. E. Knuth and O. Patashnik, Concrete Mathematics, Addison Wesley, 1993.

\bibitem{GyPsurv}
{K.~Gy{\H{o}}ry} and {{\'A}.~Pint{\'e}r}, On the equation $1^k + 2^k + \cdots + x^k = y^n$, {Publ. Math. Debrecen\/}, {62} (2003), 403--414.

\bibitem{Hirschhorn}
M. D. Hirschhorn, Evaluating $\sum_{n=1}^N(a+nd)^p$, Math. Gaz. {90} (2006), 114--116.

\bibitem{How}
{F.~T.~Howard}, Sums of powers of integers via generating functions, {Fibonacci Quart.\/}, {34} (1996), 244--256.

\bibitem{Mezo}
I. Mez\H{o}, A new formula for the Bernoulli polynomials, Result. Math. 58(3) (2010), 329--335.

\bibitem{Rahmani}
M. Rahmani, Some results on Whitney numbers of Dowling lattices,  Arab J. Math. Sci.
(2013), \\http://dx.doi.org/10.1016/j.ajmsc.2013.02.002.

\bibitem{R}
{Cs.~Rakaczki}, On the Diophantine equation $S_m (x) = g(y)$, { Publ. Math. Debrecen\/}, {65} (2004), 439--460.

\bibitem{R2}
{Cs.~Rakaczki}, On some generalizations of the Diophantine equation $s(1^k + 2^k + \cdots + x^k) + r = dy^n$, {Acta Arith. 151\/} (2012), 201-–216.

\end{thebibliography}
\end{document}